\documentclass{article}

\usepackage{graphicx}
\usepackage{dcolumn}
\usepackage{bm}

\usepackage{latexsym,amsthm}
\usepackage{amsfonts}

\newif\iffigures
\figurestrue

\newtheorem{prop}{Proposition}

\newtheorem{cor}{Corollary}

 \advance \textheight by \topskip
\begin{document}

\title{Break-up of resonant invariant circles in perturbations of the geodesic circular billiard on surfaces of constant curvature}
\author{Luciano Coutinho dos Santos \footnote{Departamento de F\'isica e Matem\'atica, CEFET-MG, Belo Horizonte, Brazil, email:astrofisico2@cefetmg.br} and S\^onia Pinto-de-Carvalho\footnote {Departamento de Matem\'atica, Universidade Federal de Minas Gerais, Brazil, email: sonia@mat.ufmg.br}}
\date{}

\maketitle
\begin{abstract}

We study  the non-persistence of horizontal invariant circles for geodesically convex perturbations of the geodesic circular billiard on surfaces of constant curvature and show that the result obtained by Ram\'irez-Ros  for the planar case \cite{rafa}, remains true for billiards  on surfaces with constant curvature.

\end{abstract}

\section{Introduction}
The plane billiard problem, defined by Birkhoff \cite{birk} in the beginning of the XX century, can be easily extended to bounded regions on surfaces of constant curvature. It is defined as the free motion of a point particle along a geodesic line in the bounded region, being reflected elastically at the impacts with the boundary. 

The motion is then completly determined  by the point of  reflection at
the boundary and the direction of movement immediately after each reflection. 
A parameter $\theta$, which locates the point of  reflection at the boundary curve, and the
angle $\psi$ between the direction of motion and the tangent to the boundary
at the reflection point, may be used to describe the system.

The billiard model defines then a map $T$ which to each $(\theta_0, \psi_0)$ in the annulus 
${\cal A} = [0,M)\times (0,\pi)$, representing the pair impact coordinate and direction of 
motion, associates the next impact and direction:
$$
\begin{array}{cccc}
 T :& {\cal A} & \to &{\cal A} \cr
    & (\theta_0, \psi_0) & \longmapsto & (\theta_1, \psi_1)
\end{array}
$$

 An oval $\Gamma$ in $S$, given in geodesic polar coordinates by $\rho=\rho(\theta)$, is a closed, simple, regular $C^q$-curve, $q\geq 2$, with strictly positive geodesic curvature. In \cite {PV}, \cite{convexgeo} and  \cite{litlle} it is proved that any oval is geodesically strictly convex, meaning that  any geodesic line cuts the oval at most twice. 
  
 If the boundary curve is an oval then the billiard map $T$ is well defined. Moreover, $T$ is a $C^{q-1}$ area preserving diffeomorphism and the billiard model gives rise to a discrete two-dimensional $C^{q-1}$ area preserving dynamical system, with the Twist Property  \cite{cons}.

In this paper, we are particularly interested on billiards on the geodesic circle, given in polar coordinates by $\rho(\theta)=\rho_0$. They are the only  billiards on surfaces with constant curvature which preserves the angle of reflection (\cite{bia}, \cite{bial}) and so their cilyndrical phase-space $ [0,M)\times (0,\pi)$ is foliated by horizontal invariant circles.

\begin{figure} [h]
\begin{center}
\includegraphics[height=.2\hsize, width=.4\hsize]{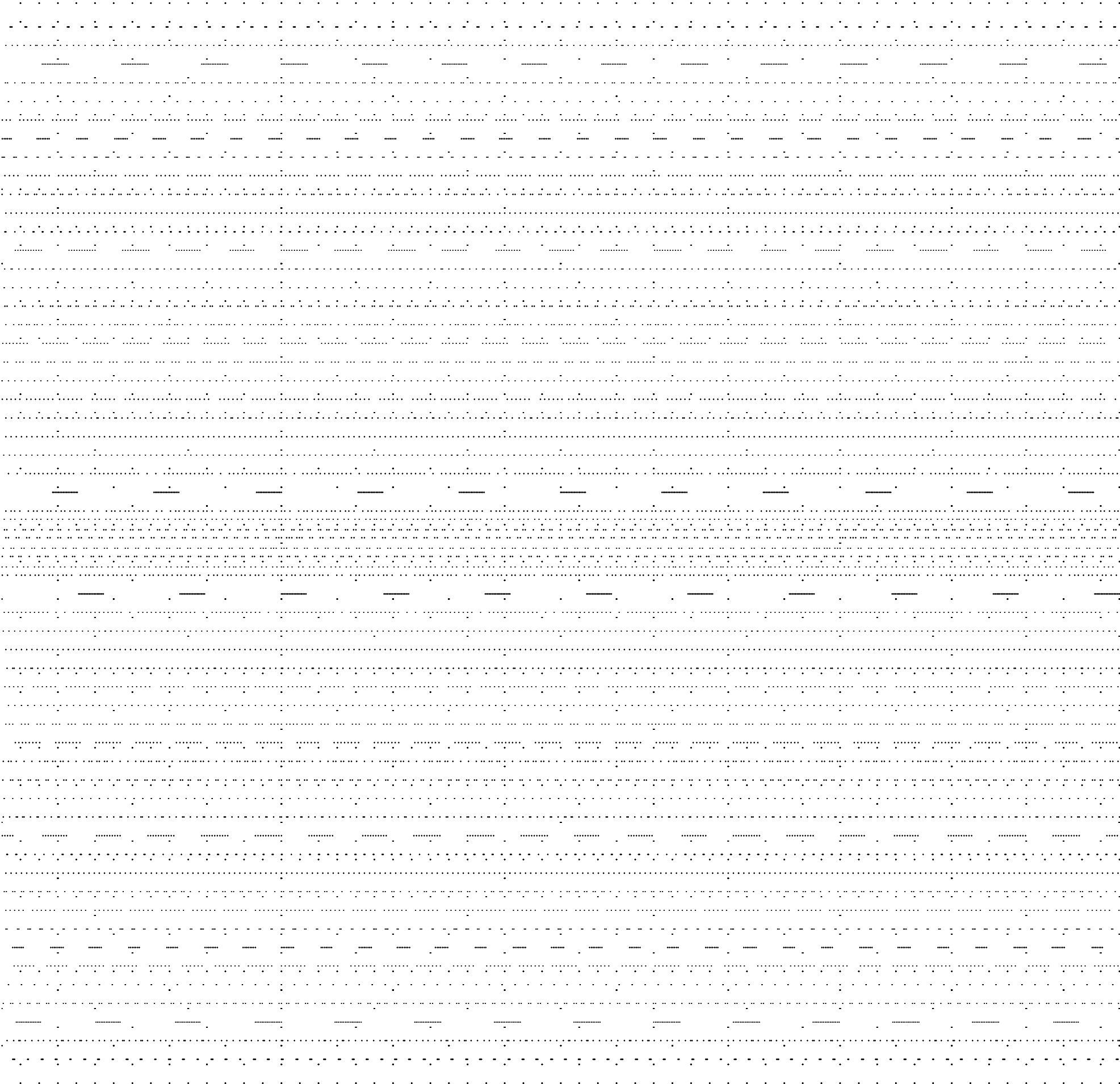}
\end{center}
\caption{Phase space of a geodesical circular billiard}
\end{figure}

Points on a same horizontal invariant circle have the same dynamical behaviour under the circular billiard map $T$. If they are periodic we say that the horizontal circle is a resonant invariant circle (RIC). 
We will show that the conditions obtained by Ram\'irez-Ros  for the planar case \cite{rafa} are still valid for the break up of a RIC under geodesically convex perturbations of the geodesic circle on a surface of constant curvature.
The main tool is the Melnikov Potential, described briefly in section \ref{sec:mel}. 

In section \ref{sec:bil} we study billiards on ovals on surfaces of constant curvature and show that the geodesic circular billiard is integrable in \ref{sec:cir}. 

Finnally , in section \ref{sec:pert} we consider perturbations of the geodesic circular billiard on the form
$$
 \Gamma_\epsilon(\theta)=\left\{\begin{array}{lcl}
 ( \rho_\epsilon(\theta) \cos\theta, \rho_\epsilon(\theta) \sin\theta,1) & \mbox{ in } & \mathbb{E}^2\\
( \sin\rho_\epsilon(\theta) \cos\theta, \sin\rho_\epsilon(\theta)\sin\theta, \cos\rho_\epsilon(\theta)) &\mbox{in}&\mathbb{S}^2_+\\
( \sinh\rho_\epsilon(\theta) \cos\theta, \sinh\rho_\epsilon(\theta)\sin\theta, \cosh\rho_\epsilon(\theta))&\mbox{in}&\mathbb{H}^2_+
\end{array}\right.
$$

with $\rho_\epsilon(\theta)=\rho_0 + \epsilon \rho_1(\theta) +  \O(\epsilon^2)$ and prove that

\begin{prop}
 If $\rho_1(\theta)=\sum\limits_{j\in\mathbb{Z}}^{} c_j \mathrm{e}^{ij\theta}$ and for $n\ge2$, there exists $j\in n\mathbb{Z}$ such that $c_j\neq 0$ then $\Upsilon_0^\frac{m}{n}$ breaks-up under $f_\epsilon$.
\end{prop}

\section{The Melnikov Potential}\label{sec:mel}

In this section we present the main steps leading to the definition of the Radial Melnikov Potential, developed by Ram\'irez-Ros in \cite{rafa} and generalized in \cite{spcrrr1}. The proofs of all the results presented here can be found in those references.

Let $\mathbb{T}=\mathbb R /2\pi\mathbb Z$, $I$ be an open interval and $\pi_1: \mathbb{T} \times I \to \mathbb{T}$ be the natural projection.
 We will use the coordinates $(x,y)$ for both $\mathbb{T} \times I$ and its cover $\mathbb R \times I$. A tilde will always denote the lift of a function or set to the universal cover.
If $g$ is a real-valued function, $\partial_i g$ denotes the derivative with respect to the $i$th variable.

Let  $f: \mathbb{T} \times I \to \mathbb{T} \times I$ be an area preserving  twist map, with generating function $g(x,x')$. Then 
$\tilde{f}(x,y)=(x',y')$ if and only if
$y  =  -\partial_1 g(x,x')$ and $y' =   \partial_2 g(x,x')$.  
Consequently, if $(x'',y'') = \tilde{f}(x',y')$, then $\partial_2 g(x,x') + \partial_1 g(x',x'') = 0.$

We study the dynamics of $f$, but it is often more convenient to work with
the lift $\tilde{f}$, so we will pass between the two without comment and,
in what follows, the lift $\tilde{f}$ remains fixed.

A closed curve $\Upsilon \subset \mathbb{T} \times I $ is said to be a
\emph{rotational invariant circle (RIC)} of $f$ when it is homotopically nontrivial and $f(\Upsilon) = \Upsilon$.
Birkhoff proved that all RICs are graphs of Lipschitz functions.
See, for instance, \cite{mei}.
Let $h: \mathbb{T} \to I$ be the Lipschitz function  such that
$\Upsilon = \mathop{\rm graph}\nolimits h := \{ (x,h(x)) : x \in \mathbb{T} \}$.
If $h$ is smooth, we say that $\Upsilon$ is a \emph{smooth RIC}.

Let $(x,y)\in \mathbb{T} \times I$ be a periodic point of the twist map $f$, and let $n$ be its least period.
Then, there exists an integer $m$ such that its lift verifies
$\tilde{f}^n(x,y)=( x + 2\pi m,y)$.
Such a periodic point is said to be of \emph{type} $(m,n)$. A RIC is said to be $(m,n)$-\emph{resonant} when all
its points are periodic of type $(m,n)$.

\begin{prop}
Let $f$ be an area-preserving twist map with a $(m,n)$-resonant smooth RIC $\Upsilon = \mathop{\rm graph}\nolimits h$ and 
 $f_\epsilon = f + \O(\epsilon)$ an area-preserving twist perturbation of $f$. 
There exist two smooth functions $h_\epsilon,h^\ast_\epsilon:\mathbb{T} \to Y$ defined for
$\epsilon \in (-\epsilon_0,\epsilon_0)$, $\epsilon_0 > 0$, such that:
\begin{enumerate}
\item
$h_\epsilon(x) = h(x) + \O(\epsilon)$ and $h^\ast_\epsilon(x) = h(x) + \O(\epsilon)$, uniformly in $x \in \mathbb{T}$; and
\item
$f_\epsilon^n\big(x,h_\epsilon(x)\big)=\big(x,h^\ast_\epsilon(x)\big)$, for all $x \in\mathbb{T}$.
\end{enumerate}
\end{prop}

We say that a $(m,n)$-resonant smooth RIC $\Upsilon$ of a twist map $f$
\emph{persists} under a perturbation $f_\epsilon = f + \O(\epsilon)$ whenever the perturbed map has
a $(m,n)$-resonant RIC $\Upsilon_\epsilon$ for any small enough $\epsilon$ such that
$\Upsilon_\epsilon = \Upsilon + \O(\epsilon)$.
Then, the resonant RIC $\Upsilon$ persists under the perturbation $f_\epsilon$
if and only if $\Upsilon_\epsilon = \Upsilon^\ast_\epsilon$.
Therefore, it is rather useful to quantify the separation
between the graphs $\Upsilon_\epsilon$ and $\Upsilon^\ast_\epsilon$,  given by
\begin{prop}
The separation between the graphs $\Upsilon_\epsilon$ and $\Upsilon^\ast_\epsilon$,  given by
$$h^\ast_\epsilon(x)- h_\epsilon(x) = L'_\epsilon(x)$$ where $L_\epsilon : \mathbb{T}\to \mathbb R$ is a function whose lift is
\begin{equation}\label{eq:CompletePotential}
\tilde{L}_\epsilon(x)=\sum_{j=0}^{n-1} g_\epsilon( \bar{x}_j(x;\epsilon),\bar{x}_{j+1}(x;\epsilon)),\qquad
\bar{x}_j(x;\epsilon)=\tilde{\pi}_1 \tilde{f}^j_\epsilon\big(x,\tilde{h}_\epsilon(x)), 
\end{equation}
and $g_\epsilon$ be  the generating function of $f_\epsilon$.
\end{prop}

It follows then that 
\begin{cor} The resonant RIC $\Upsilon$ persists under the perturbation $f_\epsilon$
if and only if $L'_\epsilon(x) \equiv 0$.\end{cor}

We shall say that $L_\epsilon : \mathbb{T} \to \mathbb R$ is the \emph{subharmonic potential} of the resonant RIC $\Upsilon$
under the twist perturbation $f_\epsilon$. It is rather natural to extract information from the low-order terms
of its expansion $L_\epsilon(x) = L_0(x) + \epsilon L_1(x) + \O(\epsilon^2)$.
This is the main idea behind any Melnikov approach to a perturbative problem.
The zero-order term $L_0(x)$ is constant (and so useless), since
$L'_0(x) = h^\ast_0(x) - h_0(x) = h(x) - h(x) \equiv 0$.
We shall say that the first-order term $L_1(x)$ is the \emph{subharmonic Melnikov potential} of the resonant RIC $\Upsilon$
under the twist perturbation $f_\epsilon$.
The proposition below provides a closed formula for its computation.

\begin{prop}\label{pro:MelnikovPotential}
If $g_\epsilon = g + \epsilon g_1 + \O(\epsilon^2)$, then the lift of $L_1(x)$ is
\[
\tilde{L}_1(x) = \sum_{j=0}^{n-1} g_1(x_j,x_{j+1}),\qquad
x_j = \tilde{\pi}_1 \tilde{f}^j(x,\tilde{h}(x)).
\]
\end{prop}

The following corollary displays the most important property of the 
subharmonic Melnikov potential in relation with the goals of this paper.

\begin{cor}\label{cor:PersistenceMelnikov}
If $L_1(x)$ is not constant,
then the resonant RIC $\Upsilon$ does not persist under the
perturbation $f_\epsilon$.
\end{cor}

\section{Billiards on ovals on surfaces of constant curvature}\label{sec:bil}

Let $S$ be a surface of constant curvature, $\Gamma\subset S$ be a closed curve and $\Omega$ be the region enclosed by $\Gamma$.
Analogously to the planar case, we can define the billiard on $\Gamma$ as the free motion of a point particle inside $\Omega$, reflecting elastically at the impacts with $\Gamma$.
Since the motion is free, the particle moves along a geodesic line of $S$ while staying inside $\Omega$ and reflects, making equal angles with the tangent at the impacts with $\Gamma$. The trajectory of the particle is a geodesic polygonal line, with vertices at the impact points.

For the study of  billiards, we will only be interested in the behaviour of the geodesics and the measure of angles. We can then take as model of surface of constant curvature $S$, one of the three surfaces:
\begin{itemize}
\item the Euclidean plane $\mathbb E^2$, given in $\mathbb R^3$ by ${\cal X}(\rho,\theta)=(\rho\cos\theta,\rho\sin\theta,1), \rho\geq0, 0\leq\theta<2\pi$;
\item an open hemisphere of the unit sphere $\S$, given in $\mathbb R^3$ by ${\cal X}(\rho,\theta)=
(\sin\rho\cos\theta,\sin\rho\sin\theta,\cos\rho), 0<\rho< \pi, 0\leq\theta<2\pi$;
\item the upper sheet of the hyperbolic plane $\mathbb{H}^2_+$, given in $\mathbb R^{2,1}$  by ${\cal X}(\rho,\theta)=
(\sinh\rho\cos\theta,\sinh\rho\sin\theta,\cosh\rho), \rho\geq0, 0\leq\theta<2\pi$.
\end{itemize}

The geodesics on $S$ are the intersections of the surface with the planes passing by the origin.  $S$ is geodesically convex and the distance between two points $X$ and $Y$ on $S$ is measured by
$$d_S(X,Y)=\, \left\{\begin{array}{ccc}
    \sqrt{<X-Y,X-Y>}\ &  \mbox{if}& X,Y\in \mathbb{E}\\
   \arccos(<X,Y>)&  \mbox{if}& X,Y\in\S\\
   \mathrm{arccosh}(-<<X,Y>>) & \mbox{if}& X,Y\in \mathbb{H}^2_+
\end{array}\right. $$

where $<,>$ is the usual inner product on $\mathbb R^3$ and $<<,>>$ is the inner product on $\mathbb R^{2,1}$.

An oval $\Gamma$ in $S$, given by $\rho=\rho(\theta)$, is a closed, simple, regular $C^q$-curve, $q\geq 2$, with strictly positive geodesic curvature. In \cite {PV}, \cite{convexgeo} and  \cite{litlle} it is proved that any oval is geodesically strictly convex, meaning that  any geodesic line cuts the oval at most twice. 

If $\Gamma$ is an oval, as $\Omega$ is a bounded subset of a geodesically convex surface, with strictly geodesically convex boundary, the billiard motion on $\Omega$ is completely determined  by the impact point and the direction of movement immediately after each reflection. 
Therefore it is determined by the points of impact at $\Gamma$, given by the polar angle $\theta$,  
and the direction of movement after each reflection, given by an angle $\psi$.

This problem defines the  Billiard Map
\begin{center}
\parbox{.4\hsize}{$
\begin{array}{ccc}
f:\mathbb{T}\times (0,\pi)&\longrightarrow &\mathbb{T}\times (0,\pi)\\
(\theta_0,\psi_0)&\longmapsto &(\theta_1,\psi_1)
\end{array}
$}\parbox{.3\hsize}{
\includegraphics[width=\hsize]{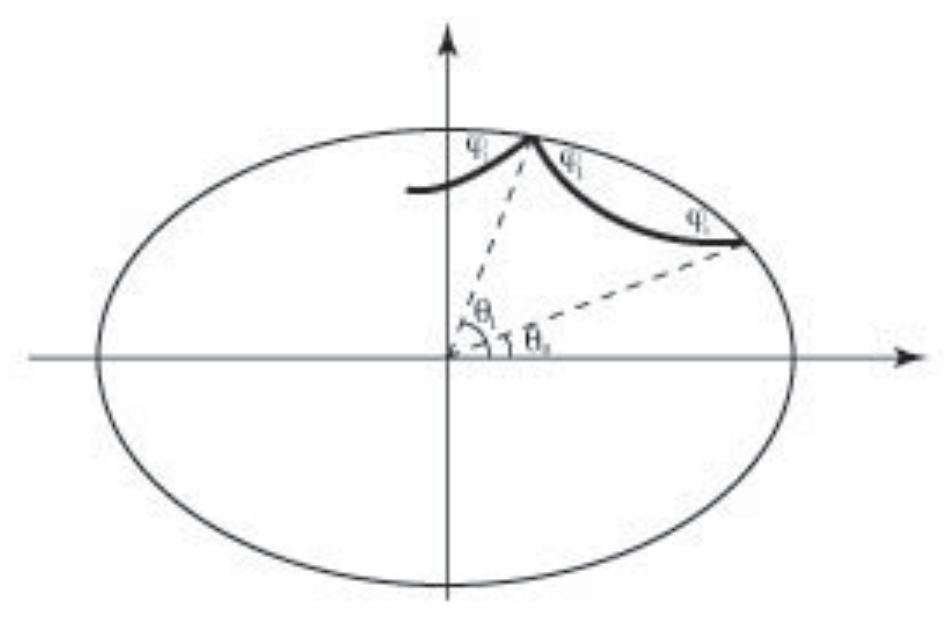}}
\end{center}

\begin{prop}
 If $\Gamma$ is a $C^q$-oval, $q\geq 2$ then the billiard map $f$ is a $C^{q-1}$-diffeomorphism, preserving the measure $d\mu = \|d\Gamma/d\theta\| \sin\psi d\psi d\theta$. It  is also a Twist map with generating function $g(\theta,\theta')=-d_S(\Gamma(\theta),\Gamma(\theta'))$.
\end{prop}
\proof In \cite{cons} we have showed that this result is true considering $\Gamma$ parameterized by the arclength parameter $s$. Since $s(\theta)=\int_0^\theta\, ||\Gamma'(\theta)||\,d\theta$, the result follows.

\section{The geodesic circular billiard is integrable}\label{sec:cir}

\parbox{.7\hsize} {
A geodesic circle $\Gamma_0$ is given by $\rho\equiv\rho_0$.  By the law of cosines on $S$, we get that the associated  billiard map is
$$f_0(\theta_0,\psi_0)= (\theta_0 + \alpha(\psi_0),\  \psi_0)$$
where
$$\alpha(\psi_0)= \left\{
\begin{array}{ll}
2\psi_0 &\mbox{in} \ \ \mathbb{E} \\
\arccos(\frac{\cos^2 \rho_0-\tan^2\psi_0}{\sec^2\psi_0-\sin^2\rho_0})&\mbox{in} \ \ \mathbb S^2_+\\
 \arccos(\frac{\cosh^2 \rho_0-\tan^2\psi_0}{\sec^2\psi_0+\sinh^2\rho_0})&\mbox{in} \ \ \mathbb{H}^2_+
\end{array}\right.$$}

It follows immediately that 
\begin{prop} The map $f_0$ is integrable and its phase-space is foliated by horizontal RIC's $\Upsilon_0(\psi_0)=\{(\theta,\psi_0), \theta \in [0,2\pi)\}$.
In particular, let $0<m<n, \gcd(m,n)=1$ and $\psi^\frac{m}{n}$ such that $\alpha(\psi^\frac{m}{n})=\frac{2\pi m}{n}$. Then $\Upsilon_0^\frac{m}{n}:= \{(\theta, \psi^\frac{m}{n}) , \theta \in [0,2\pi)\}$ is an $(m,n)$-resonant RIC.
\end{prop}

Its generating function is also easily calculated. Let $(\theta_1,\psi_1)=(\theta_0+\alpha(\psi_0), \psi_0)=f_0(\theta_0,\psi_0)$.  the generating function $g_0$ is given by
$$g_0(\theta_0,\theta_1)= \left\{
\begin{array}{ll}
-\alpha(\psi_0(\theta_0,\theta_1))\rho_0 &\mbox{in} \ \ \mathbb{E} \\
-\alpha(\psi_0(\theta_0,\theta_1))\sin\rho_0&\mbox{in} \ \ \mathbb{S}^2_+\\
-\alpha(\psi_0(\theta_0,\theta_1))\sinh\rho_0 &\mbox{in} \ \ \mathbb{H}^2_+
\end{array}\right.$$

\section{Perturbations of the geodesic circular billiard}\label{sec:pert}

We consider perturbations of the geodesic circle on the form
$$
 \Gamma_\epsilon(\theta)=\left\{\begin{array}{lcl}
 ( \rho_\epsilon(\theta) \cos\theta, \rho_\epsilon(\theta) \sin\theta,1) & \mbox{ in } & \mathbb{E}\\
( \sin\rho_\epsilon(\theta) \cos\theta, \sin\rho_\epsilon(\theta)\sin\theta, \cos\rho_\epsilon(\theta)) &\mbox{in}&\mathbb{S}^2_+\\
( \sinh\rho_\epsilon(\theta) \cos\theta, \sinh\rho_\epsilon(\theta)\sin\theta, \cosh\rho_\epsilon(\theta))&\mbox{in}&\mathbb{H}^2_+
\end{array}\right.
$$
where  $\rho_\epsilon(\theta)=\rho_0 + \epsilon \rho_1(\theta) +  \O(\epsilon^2)$. We consider $\epsilon$ is small enough so that the perturbed curve is still an oval. 
Remark that $\Gamma_\epsilon(\theta)= \Gamma_0(\theta) + \epsilon\Gamma_1(\theta)+\mathcal{O}(\epsilon^2)$, where $\Gamma_0$ is the geodesic circle with radius $\rho_0, \sin\rho_0$ or $\sinh\rho_0$.

The associated billiard map and generating function of the billiard map in $\Gamma_\epsilon$ can be written as  $f_\epsilon = f_0 + \O(\epsilon)$ and  $g_\epsilon = g_0 + \epsilon g_1+\mathcal{O}(\epsilon^2)$.

\begin{prop}
The radial Melnikov  potential of $f_\epsilon$ associated to the unperturbed RIC $\Upsilon_0^{\frac{m}{n}}$  is 
$$L_1(\theta)=C(\rho_0,m,n)\sum\limits_{j=0}^{n-1}\rho_1(\theta+\frac{2\pi m j}{n})\ \ \mbox{ with}\ \ 
C(\rho_0,m,n)= \left\{\begin{array}{lcr}
   \frac{-4 \rho_0 \sin^2\frac{\pi\, m}{n}}{l_0} & \mbox{ in }& \mathbb{E}^2\\
 \frac{2\sin(2\rho_0)\sin^2\frac{\pi m}{n}}{\sin l_0} &\mbox{in}&\mathbb S^2_+\\
  \frac{-2\sinh(2\rho_0)\sin^2\frac{\pi m}{n}}{\sinh(l_0)} &\mbox{in}&\mathbb H^2
 \end{array}\right.$$
where $l_0$ is the geodesic distance between two consecutive impacts on the nonperturbed circular billiard (with radius $\rho_0, \sin\rho_0$ or $\sinh\rho_0$, depending on $S$).
\end{prop}

\proof The case $\mathbb{E}^2$ is proved in  $\cite{rafa}$. 

Case $\mathbb{S}^2_+$: $\Gamma_\epsilon(\theta)= \Gamma_0(\theta) + \epsilon\Gamma_1(\theta)+\mathcal{O}(\epsilon^2)$, where $\Gamma_0$ is the geodesic circle with radius $\sin\rho_0$. 
\begin{eqnarray*}
\cos g_\epsilon(\theta,\bar{\theta}) &=&<\Gamma_\epsilon(\theta),\Gamma_\epsilon(\bar{\theta})>\nonumber\\
                                    &=&\sin\rho_\epsilon(\theta)\sin\rho_\epsilon(\bar{\theta})\cos\theta\cos\bar{\theta}+
                                     \sin\rho_\epsilon(\theta)\sin\rho_\epsilon(\bar{\theta})\sin\theta\sin\bar{\theta}+\cos\rho_\epsilon(\theta)
                                     \cos\rho_\epsilon(\bar{\theta})\nonumber\\
                                     &=& \sin\rho_\epsilon(\theta)\sin\rho_\epsilon(\bar{\theta})\cos(\theta-\bar{\theta})+\cos\rho_\epsilon(\theta)  \cos\rho_\epsilon(\bar{\theta}).
\end{eqnarray*}
As $\rho_\epsilon(\theta)=\rho_0 + \epsilon \rho_1(\theta) +  \O(\epsilon^2)$,
\begin{eqnarray*}
\sin\rho_\epsilon({\theta})\sin\rho_\epsilon(\bar{\theta})&=&\sin^2\rho_0 + \epsilon \sin\rho_0\cos\rho_0\,(\rho_1({\theta})+\rho_1(\bar{\theta}))
+ \mathcal{O}(\epsilon^2)\\
\cos\rho_\epsilon(\theta) \cos\rho_\epsilon(\bar{\theta}) &=& \cos^2\rho_0 - \epsilon \cos\rho_0 \sin\rho_0\,(\rho_1(\theta) + \rho_1(\bar{\theta})) + \mathcal{O}(\epsilon^2)
\end{eqnarray*}

We also have that 
  $f_\epsilon(\theta,h_\epsilon(\theta))= f_0(\theta,h_0(\theta))+\O(\epsilon)$  and $\Pi_1(f_0(\theta,h_0(\theta)))= \theta+\frac{2\pi m}{n}$. Then 
$$
\theta_j(\theta)=\Pi_1(f_\epsilon^j(\theta,h_\epsilon(\theta)))= \theta+\frac{2\pi j m}{n} +\epsilon \theta_1^j(\theta)+\mathcal{O}(\epsilon^2)
$$
Then
$$\theta_{j+1}-\theta_j=\frac{2\pi m}{n}+\epsilon \Theta^j +\mathcal{O}(\epsilon^2)\ \ \mbox{where}\ \ \Theta^j=\theta_1^{j+1}(\theta)-\theta_1^j(\theta).$$
and
\begin{eqnarray*}
\cos(\theta_{j+1}-\theta_{j})\,&=&\, \cos (\frac{2\pi m}{n} + \epsilon\Theta^j + \mathcal{O}(\epsilon^2))\nonumber\\
                           \,&=&\, \cos\frac{2\pi m}{n}\,+\,\epsilon\,\sin\frac{2\pi m}{n}\, \Theta^j +\mathcal{O}(\epsilon^2)
\end{eqnarray*}

We get then that
$$
\cos g_\epsilon(\theta_{j+1},\theta_{j})= \cos g_0 + \epsilon \beta+ \mathcal{O}(\epsilon^2)
$$
where 
$$\beta:= \Bigg[ \sin\frac{2\pi m}{n} \sin^2\rho_0\Theta^j - \sin2\rho_0 \big(\rho_1(\theta_j) + \rho_1({\theta_{j-1}})\big)\sin^2\frac{\pi m}{n} \Bigg]$$
and $\cos g_0=\sin^2\rho_0\,\cos\frac{2\pi m}{n}+\cos^2\rho_0$ is the distance between two consecutive impacts at the geodesic circle with radius $\sin\rho_0$.

Continuing we have
\begin{eqnarray*}
g_\epsilon(\theta_{j+1},\theta_{j})&=&\arccos(\cos g_0 + \epsilon \beta)+ \mathcal{O}(\epsilon^2)\\
                                    &=&g_0 + \epsilon \beta\arccos'(\cos g_0)+\mathcal{O}(\epsilon^2)
\end{eqnarray*}

and finally we get that
\begin{eqnarray*}
L_\epsilon(\theta)&=& \sum\limits_{j=0}^{n-1} g_\epsilon(\theta_{j+1},\theta_{j})\\
                  &=& \sum\limits_{j=0}^{n-1} (g_0 +\epsilon\beta\arccos'(\cos g_0)+\mathcal{O}(\epsilon^2))\\
                  &=& ng_0 + \epsilon\frac{1}{\sin g_0}\,\sum\limits_{j=0}^{n-1}\beta\, +\mathcal{O}(\epsilon^2)\\
                  &=& ng_0 +\epsilon\frac{1}{\sin g_0}\,\sum\limits_{j=0}^{n-1}\Bigg[ \sin\frac{2\pi m}{n} \sin^2\rho_0\Theta^j - \sin2\rho_0(\rho_1(\theta_j) + \rho_1({\theta_{j+1}}))\sin^2\frac{\pi m}{n} \Bigg] + \mathcal{O}(\epsilon^2)\\
                  &=& ng_0 -\epsilon\frac{\sin2\rho_0\sin^2\frac{\pi m}{n} }{\sin g_0}\sum\limits_{j=0}^{n-1}[\rho_1(\theta_j) + \rho_1({\theta_{j+1}})] + \mathcal{O}(\epsilon^2)\\
                  &=&ng_0 -\epsilon\frac{2\sin2\rho_0\sin^2\frac{\pi m}{n} }{\sin g_0}\sum\limits_{j=0}^{n-1}\rho_1(\theta_j)  + \mathcal{O}(\epsilon^2)\\
                   &=&ng_0 -\epsilon\frac{2\sin2\rho_0\sin^2\frac{\pi m}{n} }{\sin g_0}\sum\limits_{j=0}^{n-1}\rho_1(\theta+\frac{2\pi m j}{n})  + \mathcal{O}(\epsilon^2)
\end{eqnarray*}

As $l_0=-g_0$ it follows that 
$$L_1(\theta)=\frac{2\sin2\rho_0\sin^2\frac{\pi m}{n} }{\sin l_0}\sum\limits_{j=0}^{n-1}\rho_1(\theta+\frac{2\pi m j}{n}) $$

Case $\mathbb{H}^2_+$ is analogous, just taking $\sinh$, $\cosh$ and the geodesic distance in $\mathbb R^{2,1}$.

\begin{prop}
If $\rho_1(\theta)=\sum\limits_{j\in\mathbb{Z}}^{} c_j \mathrm{e}^{ij\theta}$ and for $n\ge2$, there exists $j\in n\mathbb{Z}$ such that $c_j\neq 0$ then $\Upsilon_0^\frac{m}{n}$ breaks-up under $f_\epsilon$.
\end{prop}
\proof
It suffices to remark that $\mathrm{e}^{ij(\theta+\frac{2\pi m}{n})}= n\mathrm{e^{ij\theta}}$ when $j\in n\mathbb Z$ and vanishes otherwise. So $L_1(\theta)=C(\rho_0,m,n)\sum\limits_{j\in n\mathbb{Z}}^{} c_j \mathrm{e}^{ij\theta}$. As $C(\rho_0,m,m)\neq 0$ for $m<n$, the result follows.

\vskip.5cm
Acknowledgments: The authors thank  Brazilian agencies CAPES and FAPEMIG.

\end{document}